\documentclass[a4paper, 11pt]{article}

\usepackage[T1]{fontenc}
\usepackage{amsmath}
\usepackage{amssymb}

\usepackage{amsfonts}
\usepackage{amsthm}
\usepackage{indentfirst}
\usepackage{psfrag}
\usepackage{amscd,stmaryrd,latexsym}
\usepackage[pdftex]{graphicx}
\newtheorem{thm}{Theorem}[section]

\newtheorem{Th}{Theorem}

\newtheorem*{thm*}{Theorem}

\theoremstyle{remark}

\theoremstyle{definition}

\theoremstyle{remark}

\newcommand{\be}{\begin{equation}}
\newcommand{\ee}{\end{equation}}
\newcommand{\R}{\mathbb{R}}

\newcommand{\N}{\mathbb{N}}

\def\eps{\mathop{\varepsilon}}
\def\Mc{\mathop{\mathcal{M}_{}}}

\def\Uc{\mathop{\mathcal{U}_{}}}

\def\Om{\Omega}
\def\om{\omega}

\begin{document}

\title{\textbf{Semi-calibrated $2$-currents are pseudo holomorphic, with applications.}}

\author{\textit{Costante Bellettini}\footnote{\textit{Address 1}: Fine Hall, Washington Road
Princeton NJ 08544 USA \newline \textit{E-mail address 1}: cbellett@math.princeton.edu \newline \textit{Address 2}: Einstein Drive, Princeton NJ 08540 USA \newline \textit{E-mail address 2}: costante@math.ias.edu}  \\ \\Princeton University and IAS\footnote{Research partially supported by the Giorgio and Elena Petronio Fellowship and the National Science Foundation under agreement No. DMS-1128155.}  }
\date{}
\maketitle

\textbf{Abstract}: \textit{We show that any semi-calibration of degree $2$ is locally induced by a smooth almost complex structure. We provide some applications of this result in the regularity theory for semi-calibrated $2$-currents.}

\medskip

\medskip

Calibrated geometry made its appearance, in full generality, in the landmark paper \cite{HL}. A connection with calculus of variations was immediate, since calibrated integral currents are mass minimizers in their homology class. Further interest in the topic came later on, with highly remarkable works that relate the regularity properties of calibrated currents to important geometric issues, such as theory of invariants and gauge theory (see e.g. \cite{DoT}, \cite{Ta} and \cite{Ti}). In recent years the topic has been enlarged by looking, more generally, at the weaker notion of semi-calibration (this terminology was introduced in \cite{PR}). A semi-calibration generalizes the notion of a calibration by dropping the closedness assumption. Semi-calibrated currents are no longer mass minimizers, but they turn out to be almost minimizers, in a suitable sense. The regularity theory for semi-calibrated currents is however generally expected to be better than that for arbitrary almost minimizers. For geometric applications, see  \cite{Ti}, \cite{DS}, this weaker notion yields a grater flexibility to deform the equations, which could be important for example to obtain transversality. Moreover (see \cite{HL} II.5) whenever $\phi$ is a parallel calibration of degree $k$ on $\R^n$, the contraction of $\phi$ with the Euler vector field defines a semi-calibration $\psi$ of degree $k-1$ on $\mathbb{S}^{n-1}$ and currents in the sphere semi-calibrated by $\psi$ are exactly the links of $\phi$-calibrated cones. This shows how the regularity theory for semi-calibrations has a direct impact on the regularity theory of calibrated currents themselves. Well-known examples are the Special Legendrian semi-calibration in $\mathbb{S}^{2n-1}$, that arises from Special Lagrangian geometry (see \cite{BR} for the regularity theory when $n=3$) and the nearly K\"ahler structure on $\mathbb{S}^6$, related to $G_2$ geometry and associative submanifolds. 

In this note we focus on semi-calibrations of degree $2$ and their normalized powers (the latter can be found in the final subsection). A classical example of calibration of degree two is a K\"ahler form in a complex manifold, endowed with the compatible metric (uniquely induced by the K\"ahler form and the almost complex structure). The analogue well-known semi-calibration is an almost K\"ahler form in an almost Hermitian manifold. In this case semi-calibrated integral $2$-currents turn out to be exactly the pseudo holomorphic ones, i. e. those integral $2$-currents such that the (approximate) oriented tangent planes, almost everywhere well-defined, are invariant under the action of the almost complex structure and are positively oriented (with respect to the orientation induced by the almost complex structure). This is a consequence Wirtinger's inequality. Interesting semi-calibrations of degree two also arise in the study of ruled and $2$-ruled calibrated submanifolds, e.g. \cite{Br} 3.7 and \cite{Fox}. In this work we will show that, given an arbitrary semi-calibration $\phi$ of degree $2$ on a Riemannian manifold, we can locally construct an almost complex structure on the Riemannian manifold\footnote{In the case that the ambient manifold is odd-dimensional, then we consider the cartesian product with $\R$.} that turns any $\phi$-semi-calibrated $2$-current into a pseudo holomorphic one. This can be viewed as a local classification result for semi-calibrations of degree $2$ in the spirit of the classification theory of parallel calibrations described in \cite{Joyce}, Chapter 4. The classification of parallel $2$-forms of unit comass was obtained in \cite{HL}.

In Section \ref{main} we present our main result, then in Section \ref{applications} we give a few examples of applications concerning regularity theory.

\section{How to turn semi-calibrated $2$-currents into pseudo holomorphic ones.}
\label{main}

Given a differential form $\phi$ of degree $k$ on a Riemannian manifold $\Mc$ endowed with a metric $g$, we can compare the action of $\phi$ on oriented $k$-dimensional planes with the action of the $k$-dimensional volume form. The comass $\|\phi\|^*$ of $\phi$ (see \cite{HL}) gives a quantitative meaning to this comparison, as follows. Denoting by $T_x \Mc$ the tangent space to $\Mc$ at a point $x$, recall that unit simple $k$-vectors in $T_x \Mc$ are naturally identified with oriented $k$-dimensional planes in $T_x \Mc$. Then 

\[||\phi||^*:=  \sup \{\langle \phi_x, \xi_x \rangle: x \in \Mc,\, \xi_x \text{ is a unit simple $k$-vector in } T_x \Mc\}.\] 

Forms of unit comass are called semi-calibrations (see \cite{PR}), a generalization of the notion of calibration introduced in \cite{HL} (calibrations require the closedness of the form as extra-condition). Semi-calibrated submanifolds, or more generally semi-calibrated integral currents, are defined by the condition that (almost all of) their oriented tangent planes belong to the subfamily of unit simple $k$-vectors on which $\phi$ gives its maximum action $1$, in other words $\phi$ and the $k$-dimensional volume have exactly the same action. Calibrated currents are easily seen to be mass minimizers in their homology class (see \cite{HL}). The same proof shows that semi-calibrated currents are instead almost minimizers in their homology class (or $\lambda$-minimizers, e.g. in the sense of \cite{DS}).

We show the following local result for the $2$-dimensional case.

\begin{thm}
\label{thm:main}
Given a Riemannian manifold $\Mc$ with metric $g$, let $\om$ be an arbitrary semi-calibration of degree $2$ in $(\Mc, g)$. For any point $x_0 \in \Mc$ there exist a neighbourhood $\Uc$ of $x$ in $\Mc$, a non-degenerate differential form $\Om$ of degree $2$ and a Riemannian metric $g_J$ on $\Uc$ [in the case that $\Mc$ is odd-dimensional then $\Om$ and $g_J$ are defined on $\Uc \times I$ for some interval $I \subset \R$ containing $0$] such that

\begin{itemize}

 \item $\Om$ is a semi-calibration in $(\Uc, g_J)$ [in $(\Uc \times I, g_J)$ in the case that $\Mc$ is odd-dimensional];

 \item $\Om$ and $g_J$ uniquely define an almost complex structure $J$ on $\Uc$ [on $\Uc \times I$ when $\Mc$ is odd-dimensional] such that $g_J(\cdot, \cdot) = \Om(\cdot, J \cdot)$;

 \item any $2$-plane in $(T_x \Mc, g)$ (where  $x$ is an arbitrary point on $\Uc$) that is calibrated by $\om(x)$ is also calibrated by $\Om(x)$ in $(T_x \Mc, g_J)$ [it is calibrated by $\Om(x)$ in $(T_{(x,0)} (\Mc \times \R), g_J)$ when $\Mc$ is odd-dimensional].  In particular any $2$-plane in $(T_x \Mc, g)$ (where  $x$ is an arbitrary point on $\Uc$) that is calibrated by $\om(x)$ is $J(x)$-holomorphic.

\end{itemize}
\end{thm}

In other words, by changing the semi-calibration and the metric we possibly increase the family of calibrated planes but surely preserve those that were initially calibrated. Before going into the core of the proof, let us ana\-lyse some important eigenvalue properties.

\medskip

\textbf{The endomorphism $A$ associated to $\om$ and $g$}.
We are going to define a linear map $A=A(x) :T_x \Mc \to T_x \Mc$ and study its eigenvalues in relation to the calibration properties of $\om(x)$. Set

$$\om(v,w) = g(Av, w) \text{ for any } v,w \in T_x \Mc .$$

Equivalently $Av$ is defined to be the (controvariant) vector whose $g$-dual (covariant) version is the covector $\om(v, \cdot)$. Then $A$ is a skew-adjoint operator with respect to $g$, i.e. $g(Av, w)=-g(v,Aw)$. The $g$-semipositive definite endomorphism $-A^2 = A^* A$ is thus diagonalizable. Observe that if $v$ is an eigenvector (relative to a non-zero eigenvalue) of $-A^2$ then $Av$ is an eigenvector for the same eigenvalue.

\medskip

The comass of $\om$ with respect to $g$ is $1$: taken an arbitrary vector $v \notin \text{Ker}A$ this yields, by comparing the actions of $\om$ and of the $g$-area-form on the $2$-plane $v \wedge Av$, the inequality

$$g(v,v) g(Av, Av) - g(v,Av)^2 \geq \om(v,Av)^2.$$

By definition of $A$ any vector $v$ is $g$-orthogonal to $Av$ and thus $$g(v,v) g(Av, Av)  \geq g(Av, Av)^2.$$ So, since $A$ is $g$-skew-adjoint, we can rewrite the inequality as

$$g(v,v) \geq g(-A^2v, v).$$

We can therefore conclude that any eigenvalue of $-A^2$ must belong to the interval $[0,1]$. The eigenspace of $A$ relative to $0$, i.e. $\text{Ker}A$, is the $g$-orthogonal complement of the maximal vector subspace on which $\om(x)$ is non-degenerate.

\medskip

We are now going to relate the calibrating properties of $\om(x)$ to the eigenvalues of $A$. First remark the following. Assume that $v \wedge w$ is a calibrated $2$-plane and let $t \in T_x \Mc$ be $g$-orthogonal to the $2$-plane $v \wedge w$; then $\om(v, t) = \om(w, t)=0$. This can be seen by noticing that, as $t$ varies among all possible vectors orthogonal to $v \wedge w$, the $2$-planes of the form $v \wedge t$ and $w \wedge t$ span the tangent space to the Grassmannian $G(2, T_x\Mc)$ at the point $v \wedge w$. From the fact that $\om$ restricted to $G(2, T_x\Mc)$ realizes its maximum at $v \wedge w$ we have that $\om(v, t) = \om(w, t)=0$, as desired. This fact is known as the ``first cousin principle'', see e.g. \cite{HL}, \cite{BH}.

With this in mind it follows, by the definition of $A$, that $g(Av, t) = g(Aw, t)=0$ for any $t \in T_x \Mc$ that is $g$-orthogonal to $v \wedge w$. Observe that if $v$ is a (non-zero) vector that belongs to a calibrated plane then $Av \neq 0$ by definition of $A$. Therefore $A$ restricts to an endomorphism of the $2$-plane $v \wedge w$. In other words $v \wedge Av$ is calibrated by $\om$ in $\left( T_x \Mc, g\right)$. This can be expressed by the equality

$$g(v,v) g(Av, Av) - g(v,Av)^2 = \om(v,Av)^2 .$$

Recalling that the $g$-adjoint of $A$ is $-A$ and that $v$ is $g$-orthogonal to $Av$ we rewrite the previous equality as

$$g(v,v) = g(-A^2v,v). $$

It follows that the only possible eigenvalue for the endomorphism $-A^2$ restricted to the $2$-plane $v \wedge w$ is $1$. So we conclude that, for an arbitrary choice of unit vector $v$ in a calibrated $2$-plane, the map $-A^2$ restricted to the calibrated $2$-plane is diagonalizable with $g$-orthonormal eigenbasis $\{v, Av\}$ and (double) eigenvalue $1$.

\begin{proof}[\textbf{proof of Theorem \ref{thm:main}}]
Without loss of generality we can work in a local chart around $x_0$. Therefore we will identify a neighbourhood of $x_0$ in $\Mc$ with a neighbourhood of the origin of $\R^n$. Moreover, by using the standard embedding of $\R^n$ in $\R^{n+1}$ in the case that $n$ is odd, we will further assume that we are given the semi-calibration $\om$ in a neighbourhood of $\R^{2n}$ endowed with a Riemannian metric $g$ and that $x_0$ is the origin.

The two-form $\om$ has rank $m$ at the origin for some $m \in \N$, $0 \leq m \leq n$; this means that $\om(0)^m \neq 0$ and $\om(0)^{m+1}=0$. From $\om(0)^m \neq 0$ it follows that $\om^m \neq 0$ in a neighbourhood of the origin, so the rank is at least $m$ in a neighbourhood, but it is not necessarily constant and could be strictly higher than $m$ somewhere. The rank at $x$ is half the number of non-zero eigenvalues of $A^2(x)$ (recall that the eigenvalues of $A^2$ all have multiplicity two). So at the origin we have that Ker$A(0)$ is $2(n-m)$-dimensional and there are $\lambda_1(0)\leq ... \leq \lambda_m(0) \in ]0,1]$ double eigenvalues of $A^2(0)$. Let $\eps := \lambda_1(0) >0$, so that $\lambda_j(0) \geq \eps$ for all $j=1, ..., m$.

Since $A(x)$ smoothly depends on $x$, if we stay in a small enough neighbourhood $\Uc$ of the origin there will be exactly 

\[m \text{ double eigenvalues } \lambda_1(x), ..., \lambda_m(x) \text{ of } -A^2(x) \]

that belong to the interval $[\frac{\eps}{2}, 1]$ and exactly $(n-m)$ double eigenvalues of $A^2(x)$ that belong to the interval $[0, \frac{\eps}{4}]$.

Consider the $2m$-dimensional vector subspace spanned by the eigenvectors of $A^2(x)$ relative to the eigenvalues $\lambda_1(x), ..., \lambda_m(x)$ and denote it by $V(x)$. The $g$-orthogonal complement $V^\bot(x)$ is spanned by the eigenvectors relative to the remaining eigenvalues, i.e. the eigenvalues that are $\leq \frac{\eps}{4}$. The smooth dependence of $A(x)$ on $x$ yields that $V(x)$ and $V^\bot(x)$ also depend smoothly on the point $x$.

We are going to change the metric $g(x)$ to a new metric $g_J(x)$ in the following way. Remark that $-A^2(x)$ is an endomorphism of $V(x)$, since $V(x)$ is spanned by eigenvectors of $-A^2(x)$. We can diagonalize the $g(x)$-positive definite map $-A^2(x)$ restricted to $V(x)$ and by taking the square roots of the entries we find a $g(x)$-positive definite $Q(x)$ on $V(x)$ such that $Q^2(x) = -A^2(x)$. Observe that $Q$ (and in the same way also $Q^{-1}$) commutes with $A$. To see this, choose an eigenbasis of $V(x)$ (for $-A(x)^2$) of the form $\{v_1, A v_1, ..., v_m, A v_m\}$, where $v_i$ and $A v_i$ are eigenvectors relative to $\lambda_i$. In this basis $Q$ and $A$ are $2 \times 2$-block diagonal matrices: for $Q$ each block is $\displaystyle  \left(\begin{array}{cc}
                                                                                                    \sqrt{\lambda_i} & 0 \\
0 & \sqrt{\lambda_i}                                                                                                                                                                                                                                                                                                                                                                                                                                                                                                                                                                                                                                                                                        
                                                                                                                                                                                                                                                                                                                                                                                                                                                                                                                                                                                                                                                                                                                                                                                           \end{array}\right)$
 and the corresponding block for $A$ is $\displaystyle \left(\begin{array}{cc}
                                                                                                    0 & -\lambda_i \\
1 & 0                                                                                                                                                                                                                                                                                                                                                                                                                                                                                                                                                                                                                                                                                        
                                                                                                                                                                                                                                                                                                                                                                                                                                                                                                                                                                                                                                                                                                                                                                                           \end{array}\right)$. So they commute. Therefore by setting $J(x):=Q^{-1}(x)A(x)$ we define an endomorphism such that $J^2=-Id$, i.e. $J(x)$ is an almost complex structure on $V(x)$. Then define a new (smooth) Riemannian metric $g_J(x)$ on $\Uc$ as follows 

\be
\begin{split}
\left(g_J(x)\right)(v, w):= \om(x)(v, J(x) w) = g(x)(A(x)v, J(x)w) \text{ for } v,w \in V(x) \\
\left(g_J(x)\right)(v, w):=0 \text{ for } v \in V, w \in V(x)^\bot \\
\left(g_J(x)\right)(v, w):= g(x)(v, w) \text{ for } v,w \in V(x)^\bot.
\end{split}
\ee

The restriction of $\om(x)$ to $V(x)$ is a non-degenerate two-form on $V(x)$ that is a calibration with respect to the metric $g_J(x)$ and the calibrated $2$-planes are exactly those that are $J(x)$-holomorphic: indeed in each $V(x)$ we have a standard symplectic form related to the metric and the almost complex structure by $g_J(v, w)=\om(v, Jw)$.

The important remark here is that any $2$-plane that is calibrated by $\om(x)$ in $(T_x \Mc g)$ must lie in $V(x)$, since as we saw above it must live in the eigenspace of $-A^2(x)$ relative to the eigenvalue $1$. The map $-A^2$ restrics to an endomorphism of any calibrated $2$-plane and it is there diagonalizable with $g$-orthonormal eigenbasis of the form $\{v, Av\}$ and (double) eigenvalue $1$. In particular $Q$ is the identity on any calibrated $2$-plane. This means that, on a calibrated $2$-plane, $J=A$ and $g_J = g$ and we can see that any $2$-plane that is calibrated by $\om(x)$ in $(T_x \Mc g)$ will be calibrated as well by the restriction of $\om(x)$ to $\left(V(x), g_J(x)\right)$.

\medskip

We observe now that by restricting the action of $\om$ to the $V(x)$ we can define a smooth two-form in $\Uc$. Precisely, denoting by $\Pi(x)$ the orthogonal projection from $T_x \Mc$ to $V(x)$, we define

$$\om_1(v,w) := \om(\Pi v, \Pi w)  = g(A \Pi v, J \Pi w)$$

and this depends smoothly on $x$ since so do $g$, $\Pi$, $A$ and $J$. The action of $\om_1$ on $V(x)$ is exactly that of $\om(x)|_{V(x)}$, while $\om_1$ is degenerate on the complementary subspace $V(x)^\bot$.  

We now choose a smooth local orthonormal frame $\{t_1(x), ..., t_{2n-2m}(x)\}$ for the subspaces $V^\bot(x)$ for $x \in \Uc$ (maybe for this we need to make the initially chosen neighbourhood smaller, but by abuse of notation we keep denoting it by $\Uc$). Denote by $\star$ the metric duality between vectors and covectors with respect to $g_J$. The form

$$\om_2(x):= \left(t_1(x)\right)^\star \wedge \left(t_2(x)\right)^\star +\left(t_3(x)\right)^\star \wedge \left(t_4(x)\right)^\star+ ... \left(t_{2n-2m-1}(x)\right)^\star \wedge \left(t_{2n-2m}(x)\right)^\star $$

is smooth and degenerate on $V(x)$. Extend further the almost complex structure $J$ from $V(x)$ to the whole of on $T_x \Mc$  by setting

$$J(t_1(x))=t_2(x), ... , J(t_{2n-2m-1}(x))=t_{2n-2m}(x).$$

The two-form $\om_2$ restricts to a standard symplectic form on $V(x)^\bot$ with compatible metric $g_J$ and almost complex structure $J$. We can now define the smooth two-form

$$\Om:= \om_1 + \om_2.$$

Since $\om_1$ restricts to a standard symplectic form on $V(x)$ with compatible metric $g_J$ and almost complex structure $J$ we can see that $\Om$ has unit comass with respect to $g_J$ (e.g. by Wirtinger's inequality) and $g_J(\cdot, \cdot) = \Om(\cdot, J \cdot)$. Moreover we already pointed out that any $2$-plane that is calibrated by $\om(x)$ in $(T_x\Mc, g)$ must lie in $V(x)$ and is therefore calibrated by $\om_1(x)$ in $(V(x), g_J)$. The restriction of $\Om(x)$ to $V(x)$ is exactly $\om_1(x)$ and thus on $(V(x), g_J)$ the $2$-planes calibrated by $\om_1$ and $\Om$ are the same. So what happens is that the $2$-planes that are calibrated by $\om(x)$ in $(T_x\Mc, g)$ are a subset of those that are calibrated by $\Om(x)$ in $(T_x\Mc, g_J)$. 

\end{proof}

\section{Applications.}
\label{applications}

The theorem in Section \ref{main} shows that, given any semi-calibrated current of dimension $2$, we can locally find a new Riemannian metric $g_J$ and a new semi-calibration $\Om$ (of unit comass with respect to the new metric) such that $g_J$ and $\Om$ uniquely define an almost complex structure $J$ which makes the $2$-dimensional current (which stays untouched) a pseudo holomorphic one. In other words we can always locally reach a situation where we work with a pseudo holomorphic current in an almost Hermitian manifold. This allows us to have more structure without changing the local regularity properties of the current. We now present a few examples of the applicability of Theorem \ref{thm:main}.

\medskip

\textbf{Study of the regularity of an arbitrary semi-calibrated integral $2$-cycle}. In \cite{RT2} T. Rivi{\`e}re and G. Tian study the regularity of pseudo holomorphic cycles under a ``locally-taming'' assumption, i.e. they assume that there exists locally a closed non-degenerate two-form that is compatible with the almost complex structure $J$. They prove that such a cycle can only have isolated singularities. It would be very important to improve (and possibly shorten the proof of) their result by dropping the ``locally-taming'' assumption, as suggested by the authors themselves in the introduction of their paper. 

The case of a calibrated integral $2$-cycle embeds (thanks to the closedness assumption on the $2$-form) in the theory of mass-minimizing integral currents studied in \cite{Alm} and \cite{Ch}. Remark that in a series of recent works a new proof of this regularity theory is presented, see  \cite{DeLelSpa}, \cite{DeLelSpa4}, \cite{DeLelSpa5}, \cite{DeLelSpa6}. For the case of semi-calibrated integral $2$-cycles, however, the optimal regularity result is only known in some particular cases, e.g. Special Legendrian cycles (see \cite{BR}, \cite{B}). In this latter case the proof actually argues, in one of its first steps, by finding an almost complex structure that induces the semi-calibration.

By Theorem \ref{thm:main} the extension of \cite{RT2} (to the case when the almost complex structure is not locally tamed by a closed form) would actually yield the \underline{complete} regularity theory for semi-calibrated integral $2$-cycles. This would be a very important result in view of the geometric applications outlined in \cite{Ti}, Section 6. A possible first step to approach the problem when we drop the ``locally-taming'' assumption could be the technique used in \cite{B4} to prove the rate of decay of the mass ratio.

\medskip

\textbf{Simplification in \cite{PR}}. In \cite{PR} D. Pumberger and T. Rivi{\`e}re prove the following

\begin{thm*}
An arbitrary semi-calibrated integral $2$-cycle possesses at any point a unique tangent cone.
\end{thm*}

The first part of the paper is devoted to the proof of the theorem in the case that the semi-calibration $\om$ has an associated compatible almost complex structure $J$ such that any semi-calibrated integral $2$-cycle is actually $J$-pseudo holomorphic. In the last part of the paper the authors deal with the general case of a semi-calibrated integral $2$-cycle by approximating, in a suitable way, any such cycle with a pseudo holomorphic one. 

In view of our Theorem \ref{thm:main} this approximation can be completely avoided, since we can define an almost complex structure $J$ that makes the given semi-calibrated cycle $J$-pseudo holomorphic. The change of metric has no effect on the uniqueness property that we need to prove for the tangent cone.

\medskip

\textbf{Improvements in \cite{B3}}. \textbf{(i)} In the same vein, our Theorem \ref{thm:main} shows that the \textit{pseudo algebraic blow up} technique used in \cite{B2}, \cite{B3} to prove the uniqueness of tangent cones for a positive-$(1,1)$ integral cycle can actually be used to show the uniqueness of tangent cones for an arbitrary semi-calibrated $2$-cycle, thereby recovering, with different techniques, the same result as \cite{PR}. Remark that, in view of \cite{B4}, the \textit{pseudo algebraic blow up} technique also yields an explicit rate of decay for the mass-ratio of semi-calibrated two-cycles under a blow up analysis.

\textbf{(ii)} Remark that, given a semi-calibration $\om$ of degree $2$ in a Riemannian manifold $(\Mc, g)$, the $2p$-form $\frac{1}{p!}\om^p$ is still a semi-calibration. This can be seen by the following considerations, that show (among other things) that the comass of $\frac{1}{p!}\om^p$ is $1$. It is enough to do this pointwise so we fix a point $x$. Let $W(x) \subset T_x \Mc$ be the subspace where $\om(x)$ is non-degenerate. Remark that it only makes sense to consider values of $p$ that are not greater than the rank of $\om$ and that $\frac{1}{p!}\om^p$ is the zero-form on the orthogonal to the subspace $W(x)$, i.e. we can think of both $\om(x)$ and $\frac{1}{p!}\om(x)^p$ as forms on $W(x)$.  We perform the construction of the proof of Theorem \ref{thm:main} in the subspace $W(x)$, i.e. the $g$-orthogonal to $\text{Ker}A(x)$. This gives us the new metric $g_J$ so that $\om$ is a semi-calibration in $(W(x), g_J)$. Wirtinger's inequality tells us that $\frac{1}{p!}\om^p$ is a semi-calibration in $(\Mc, g_J)$. From the proof in Section 1 we get that $g \geq g_J$ on $W(x)$, since the eigenvalues of $-A(x)^2$ on $W(x)$ belong to $]0,1]$. This fact together with Theorem 6.11 of \cite{HL} yields that $\frac{1}{p!}\om^p$ is a semi-calibration also in $(W(x), g)$ and moreover that the set of calibrated $2p$-planes in $(W(x), g)$ is contained in the set of calibrated $2p$-planes in $(W(x), g_J)$; in particular the set of calibrated $2p$-planes in $(W(x), g)$ is contained the the eigenspace of $-A(x)^2$ relative to the eigenvalue $1$ and is made of $2p$-planes that are invariant under the action of $A(x)$ (recall that $A(x)=J(x)$ on this eigenspace). 

By passing from $\om$ to $\Om$ as in Section 1, and thus from $\frac{1}{p!}\om^p$ to $\frac{1}{p!}\Om^p$ we enlarge the set of $2p$-calibrated planes: indeed the $2p$-planes calibrated by $\frac{1}{p!}\Om^p$ in $(T_x \Mc, g_J)$ are exactly those that are $J(x)$-holomorphic. By what we previously observed, they thus contain in particular the $2p$-planes calibrated by $\frac{1}{p!}\om^p$ in $(W(x), g)$. 

So we have concluded that the $2p$-planes calibrated by $\frac{1}{p!}\om^p$ in $(\Mc, g)$ (i.e. the semi-calibration that we had in the beginning) are also calibrated by $\frac{1}{p!}\Om^p$ in $(T_x \Mc, g_J)$. The uniqueness of tangent cones for a $2p$-integral cycle $T$ semi-calibrated by $\frac{1}{p!} \om^p$ in $(\Mc, g)$ can be thus deduced from the uniqueness of tangent cones for $T$ when we view this cycle as being semi-calibrated by $\frac{1}{p!}\Om^p$ in $(\Mc, g_J)$ (for the latter case we can use the result of \cite{B3}). Therefore we have proved

\begin{thm}
Let $\om$ be a semi-calibration of degree $2$ in $(\Mc, g)$ and consider the semi-calibration of degree $2p$ given by $\displaystyle \frac{1}{p!}\om^p$, where $p\leq \frac{\text{dim}\Mc}{2}$ is a positive integer. Then any $(2p)$-dimensional integral cycle semi-calibrated by $\displaystyle \frac{1}{p!}\om^p$ possesses at any point a unique tangent cone.
\end{thm}

\end{document}